# THE SURVIVAL OF LARGE DIMENSIONAL THRESHOLD CONTACT PROCESSES

By Thomas Mountford[1] and Roberto H. Schonmann[2]

*École Polytechnique Fédérale and University of California at Los Angeles*

We study the threshold $\theta$ contact process on $\mathbb{Z}^d$ with infection parameter $\lambda$. We show that the critical point $\lambda_c$, defined as the threshold for survival starting from every site occupied, vanishes as $d \to \infty$. This implies that the threshold $\theta$ voter model on $\mathbb{Z}^d$ has a nondegenerate extremal invariant measure, when $d$ is large.

**1. Introduction and results.** In this paper we study the threshold $\theta$ contact process on the $d$-dimensional cubic lattice $\mathbb{Z}^d$. [See, e.g., Liggett (1985) and Liggett (1999) for background on the area of interacting particle systems, to which this process belongs.] The threshold $\theta$ contact process on $\mathbb{Z}^d$ is defined as the Markov process on $\{0,1\}^{\mathbb{Z}^d}$ with flip rates at site $x \in \mathbb{Z}^d$ at time $t \geq 0$ given by:

- 1 flips to 0 at rate 1.
- 0 flips to 1 at rate $\lambda$ in case there are at least $\theta$ nearest neighbors of $x$ in state 1 at time $t$; and at rate 0 otherwise.

The parameter $\lambda \geq 0$ is called the infection rate. The state of the process at each site at each time is called the spin at that site at that time. A spin 0 is interpreted as a vacant or healthy site, while a spin 1 is interpreted as an occupied or infected site. It is well known [see, e.g., Chapter I of Liggett (1985)] that such rates define a unique Markov process. Note that the flip rates above are attractive [see Chapter III of Liggett (1985)], a property that has many consequences, and will be used extensively in this paper.

When $\theta = 1$, the threshold contact process is easier to analyze, among other reasons because it has an additive dual process [see Chapter III of

Received June 2008; revised October 2008.
[1]Supported in part by Grant SNF 200020-115964.
[2]Supported in part by NSF Grant DMS-03-00672.

Liggett (1985) for additive duality]. (In this dual process infected sites become healthy at rate 1, and they infect simultaneously all their neighbors at rate $\lambda$.) Also, only when $\theta = 1$ can every site become infected when the process is started from a finite number of infected sites. The behavior of the process with $\theta = 1$ is therefore not expected to be qualitatively different from that of the much studied (linear) contact process [see Chapter VI of Liggett (1985) and Part I of Liggett (1999)]. In contrast, when $\theta \geq 2$, not only is the analysis of the model more challenging, but also the behavior is expected to be qualitatively distinct from that of the threshold 1 or the linear contact process. We defer to Section 5, item (2), a discussion of this point, and we refer the reader to Toom (1974), Durrett and Gray (1990), Bramson and Gray (1991), Chen (1992, 1994) and Fontes and Schonmann (2008) for background on threshold contact processes with $\theta \geq 2$ and related models. Here we will extend to every value of $\theta$ a result that is already available in the case $\theta = 1$.

Threshold contact processes can also be seen as natural generalizations of the bootstrap percolation processes [see, e.g., Adler and Lev (2003), for a review], with bootstrap percolation corresponding to the case with $\lambda = \infty$ (or in other words, with the down flips suppressed).

For coupling purposes, it is convenient to construct the process using a system of Poisson marks. For this purpose, associate to each site in $\mathbb{Z}^d$ two independent Poisson processes: one with rate 1, and one with rate $\lambda$. Mark the arrival times of the former with symbols $D$ (for "down") and those of the latter with symbols $U$ (for "up"). Make these Poisson processes independent from site to site. Use the marks now in the obvious way, to define the process: A spin 1 at site $x$ flips to 0 when it encounters a $D$ mark there; a spin 0 at site $x$ flips to 1 when it encounters an $U$ mark there and at least $\theta$ neighbors of $x$ have spin 1 at that time. The probability space on which these Poisson processes are defined will be large enough to accommodate the process started from arbitrary initial configurations. The same system of Poisson marks will be used to build various comparison processes that will be introduced through the paper.

We will denote by $(\eta^\mu_{\mathbb{Z}^d,\theta,\lambda;t})_{t\geq 0} = (\eta^\mu_t)_{t\geq 0}$ the process started from a random configuration picked according to law $\mu$ at time 0. When $\mu$ is concentrated on a single configuration $\zeta$, we will use the notation $(\eta^\zeta_t)_{t\geq 0}$ for the process. When $\mu$ is product measure with density $p$ we will use the notation $(\eta^p_t)_{t\geq 0}$. Similar shorthand notation will also be used in other places below.

The point mass on the configuration with all sites in state $i \in \{0, 1\}$ will be denoted $\delta_i$. The distribution $\delta_0$ is trivially invariant for the threshold $\theta$ contact process, when $\theta > 0$. By attractivity, $\eta^1_t \Rightarrow \bar{\nu}$, as $t \to \infty$, where $\Rightarrow$ denotes convergence in distribution, and $\bar{\nu}$, is called the upper invariant measure.

We say that the process started from the distribution $\mu$ dies out when $\eta_t^\mu \Rightarrow \delta_0$, as $t \to \infty$. When this happens for every $\mu$, we simply say that the process dies out. Attractivity implies that the process dies out precisely when $\bar{\nu} = \delta_0$. When the process does not die out, we will say that it survives.

Set $\rho_t^\mu = \mathbb{P}(\eta_t^\mu(0) = 1)$. The critical point for survival is defined as

$$\lambda_c = \sup\{\lambda : \bar{\nu} = \delta_0\} = \sup\{\lambda : \rho_t^1 \to 0 \text{ as } t \to \infty\}.$$

It is a standard, easy fact that $\lambda_c > 0$, and it is also known that for each $d \geq 2$ and $\theta \leq d$, $\lambda_c < \infty$. Durrett and Gray (1990) have an explicit proof in the case $\theta = 2$, using contour methods. But a well-known argument, that we present in Section 5, item (5), allows one to use the results of Toom (1974) and Bramson and Gray (1991) in discrete time to obtain this result in continuous time, for arbitrary $\theta \leq d$. Here we will prove the following.

THEOREM 1.1. *For any $\theta$, the threshold $\theta$ contact process has*

$$\lim_{d \to \infty} \lambda_c = 0.$$

In case $\theta = 1$ this result is not new. As pointed out in the display in page 191 in Cox and Durrett (1991), in this case it is known that $\lambda_c \leq C/d$, for an appropriate constant $C$. The argument given there comes from Holley and Liggett (1978), who used it in the case of the linear contact process.

From the analysis of a mean-field version of the threshold $\theta$ contact process, in Fontes and Schonmann (2008), it is natural to conjecture that for each $\theta$ there is a constant $\Phi_\theta \in (0, \infty)$ such that $d\lambda_c \to \Phi_\theta$, as $d \to \infty$. Our methods unfortunately are not sharp enough to prove this strengthening of Theorem 1.1, or even to obtain good estimates on how fast $\lambda_c$ vanishes, as $d \to \infty$.

Theorem 1.1 has an implication for the threshold $\theta$ voter model on $\mathbb{Z}^d$, defined by the following rates. At $x \in \mathbb{Z}^d$ at time $t \geq 0$:

- 1 flips to 0 at rate 1 in case there are at least $\theta$ nearest neighbors of $x$ in state 0 at time $t$; and at rate 0 otherwise.
- 0 flips to 1 at rate $\lambda$ in case there are at least $\theta$ nearest neighbors of $x$ in state 1 at time $t$; and at rate 0 otherwise.

When $\theta \geq 1$, the threshold $\theta$ voter model obviously has $\delta_0$ and $\delta_1$ as two of its extremal invariant measures. In case $\theta \geq 2$, it is also easy to find infinitely many other extremal invariant measures concentrated on a single configuration that is a trap for the process. On the other hand, from the work of Cox and Durrett (1991), restated as Proposition 2.11 in Part II of Liggett (1999), we know that a sufficient condition for this process to have a nondegenerate (i.e., not concentrated on a single configuration) extremal invariant measure is that the threshold $\theta$ contact process survives when $\lambda = 1$. In combination with Theorem 1.1 above, we obtain

THEOREM 1.2. *For any $\theta$, the threshold $\theta$ voter model on $\mathbb{Z}^d$ has a nondegenerate extremal invariant measure when $d$ is large.*

In the case $\theta = 1$ a much more detailed result is available, from Liggett (1994): the threshold $\theta = 1$ voter model on $\mathbb{Z}^d$ has a nondegenerate extremal invariant measure if an only if $d \geq 2$.

The proof of Theorem 1.1 will be split into Sections 2, 3 and 4. The proof will be presented in the case $\theta = 2$, for notational simplicity, and so that the reader can more easily visualize the arguments. In Section 5 we will explain why the proof applies to the other values of $\theta$. Section 5 will also contain other results and comments related to Theorem 1.1 and its proof.

**2. Approach to the proof of Theorem 1.1.** As mentioned in the Introduction, in this section and the following two, we suppose that $\theta = 2$.

We will compare the threshold 2 contact process on $\mathbb{Z}^d$ with the same processes on subsets of $\mathbb{Z}^d$, with free boundary conditions. For this purpose we will use the following notation. For each $x \in \mathbb{Z}^d$, let $\mathcal{N}_x$ be the set of its $2d$ nearest neighbors, and given $R \subset \mathbb{Z}^d$, set $\mathcal{N}_{R,x} = \mathcal{N}_x \cap R$. Given $R \subset \bar{R} \subset \mathbb{Z}^d$, and $\zeta \in \{0,1\}^{\bar{R}}$, we will denote by $\eta_{R;t}^\zeta$ the threshold 2 contact process on $R$, started from the configuration $\zeta$ restricted to $R$. In this process, at site $x \in R$, at time $t \geq 0$:

- 1 flips to 0 at rate 1.
- 0 flips to 1 at rate $\lambda$ in case there are at least 2 sites in $\mathcal{N}_{R,x}$ in state 1 at time $t$; and at rate 0 otherwise.

Define $G = \mathbb{Z}^2 \times \{0,1\}^{d-2}$. We will denote by $\lambda_c(G)$ the threshold for survival (started from all sites occupied) for the threshold 2 contact process on $G$. We will prove

PROPOSITION 2.1.
$$\lim_{d \to \infty} \lambda_c(G) = 0.$$

By attractivity, $\lambda_c(\mathbb{Z}^d) \leq \lambda_c(G)$, so this result implies the case $\theta = 2$ of Theorem 1.1. In this section we will explain how the proof of Proposition 2.1 will be reduced to the proof of two lemmas.

Fix a value of $\lambda > 0$. We want to prove that for $d$ large enough, the process $(\eta_{G;t}^1)_{t \geq 0}$ survives. For this purpose, we will partition $\mathbb{Z}^2 \times \{0,1\}^{d-2}$ into copies of the hypercube $H = \{0,1\}^d$, indexed by the elements of $\mathbb{Z}^2$. The elements of this partition will be denoted by

$$H_i = \{2i_1, 2i_1 + 1\} \times \{2i_2, 2i_2 + 1\} \times \{0,1\}^{d-2}, \qquad i = (i_1, i_2) \in \mathbb{Z}^2.$$

We will observe the process $(\eta^1_{G;t})_{t\geq 0}$ at times that are multiples of $d^2$. At each one of these times, we will decide whether each hypercube $H_i$, $i \in \mathbb{Z}^2$, is in a good configuration or not according to a technical definition provided later in this section. (For large $d$, the completely empty configuration will not be a good configuration, while the completely full configuration will be a good configuration.) This will be done in a way that will allow us to use a classic result of Toom (1974), to show that the hypercube $H = H_0$ will be in a good configuration at arbitrarily large times, implying therefore that $(\eta^1_{G;t})_{t\geq 0}$ survives.

We recall now the result from Toom (1974), that will be used. Let $X_{i,n}$, $i \in \mathbb{Z}^2$, $n \in \mathbb{Z}$, be a collection of i.i.d. random variables, with $P(X_{i,n} = 0) = \varepsilon$, $P(X_{i,n} = 1) = 1 - \varepsilon$, where $\varepsilon \in [0,1]$ is a noise parameter. Define now the following family of random transformations, $T_{\varepsilon;n}$, on $\{0,1\}^{\mathbb{Z}^2}$, indexed by $n \in \{0, 1, 2, 3, \ldots\}$. $T_{\varepsilon;0}$ is the identity and recursively, for $n = 1, 2, \ldots$,

$$(T_{\varepsilon;n}\Psi)_i = \begin{cases} X_{i,n}, & \text{if } (T_{\varepsilon;n-1}\Psi)_i = 1, \text{ or} \\ & (T_{\varepsilon;n-1}\Psi)_{i+(0,1)} = (T_{\varepsilon;n-1}\Psi)_{i+(1,0)} = 1, \\ 0, & \text{otherwise,} \end{cases}$$

for $\Psi \in \{0,1\}^{\mathbb{Z}^2}$, $i \in \mathbb{Z}^2$. Denote by $\mathbf{1}$ the element of $\{0,1\}^{\mathbb{Z}^2}$ identically 1. Toom (1974) proved the following result [see also Bramson and Gray (1991), for an alternative proof]: There exists $\varepsilon_c > 0$ such that

(2.1) $\qquad\qquad$ for $\varepsilon < \varepsilon_c, \quad \liminf_{n \to \infty} P((T_{\varepsilon;n}\mathbf{1})_0 = 1) > 0$.

Returning to the threshold 2 contact process on $G$, we now explain how good configurations are defined. For this purpose, consider the threshold 2 contact process on $H$, started from $\zeta \in \{0,1\}^G$. We want to compare this process with the system of independent spins on $H$, that starts from the same configuration, and uses the same system of up and down marks, for its flips. (So that flips down occur at rate 1 and flips up occur at rate $\lambda$, independently of the state of the other spins.) Call this independent flip process $(\pi^\zeta_{H;t})_{t\geq 0}$. Clearly, $\eta^\zeta_{H;t} \leq \pi^\zeta_{H;t}$, for all $t \geq 0$. Roughly speaking, a good configuration $\zeta$ on $H$, is one for which the processes $(\eta^\zeta_{H;t})_{0 \leq t \leq d^2}$ and $(\pi^\zeta_{H;t})_{0 \leq t \leq d^2}$ are likely to be close to each other. The technical notion of closeness that will be used is defined next.

DEFINITION (*Thin sets*). Given $R, S \subset \mathbb{Z}^d$, and $L > 0$, say that $S$ is $L$-thin in $R$, if

$$|S \cap \mathcal{N}_{R,x}| \leq L,$$

for all $x \in R$.

DEFINITION (*Good configurations*).   Given $R \subset \bar{R} \subset \mathbb{Z}^d$, and $\zeta \in \{0,1\}^{\bar{R}}$, we define the discrepancy sets by

$$\Delta_{R;t}^{\zeta} = \{x \in R : \eta_{R;t}^{\zeta}(x) \neq \pi_{R;t}^{\zeta}(x)\}$$
$$= \{x \in R : \eta_{R;t}^{\zeta}(x) = 0, \pi_{R;t}^{\zeta}(x) = 1\}, \qquad t \geq 0,$$
$$\Delta_R^{\zeta} = \bigcup_{t \in [0, d^2]} \Delta_{R;t}^{\zeta}.$$

And given also $L > 0$, and $\delta \in [0,1]$, we say that the configuration $\zeta$ is $(L, \delta)$-good on $R$, if

$$\mathbb{P}(\Delta_R^{\zeta} \text{ is } L\text{-thin on } R) \geq 1 - \delta.$$

The next lemma will be the key to providing us with good configurations. First we need one more technical definition.

DEFINITION (*Almost product distributions*).   Given $R \subset \mathbb{Z}^d$, $M > 0$ and $\alpha \in [0,1]$, say that a probability distribution $\nu$ over $\{0,1\}^R$ is an $M$-almost $\alpha$-product distribution, if it can be constructed as follows. There is a probability space on which there is a random configuration $\xi \in \{0,1\}^R$, with law given by an $\alpha$-product distribution (i.e., each site of $R$ is independently occupied in $\xi$, with probability $\alpha$). On this probability space there is also a random set $S \subset R$ that is almost surely $M$-thin on $R$. The law $\nu$ can then be obtained as the law of the random configuration defined as being identical to $\xi$ on $R \setminus S$ and identically 0 on $S$. (Note that in this definition, the set $S$ does not need to be independent of $\xi$.)

CONVENTION ON NOTATION AND TERMINOLOGY.   In the three definitions above, when $R = H$, we may omit it. For instance, "$\zeta$ is good" has the same meaning as "$\zeta$ is good on $H$," $\Delta^{\zeta} = \Delta_H^{\zeta}$, etc.

LEMMA 2.1 (Almost product distributions concentrate on good configurations).   *For any $\lambda > 0$ and any $\alpha \in (0,1]$ there is $L = L(\lambda, \alpha) = L(\lambda/(1+\lambda) \wedge \alpha)$ such that for any $M > 0$, $\delta > 0$, and $\varepsilon > 0$, there is a finite $d_0$ such that if $\nu$ has an $M$-almost $\alpha$-product distribution over $\{0,1\}^H$, then*

$$\nu(\{\zeta : \zeta \text{ is } (L, \delta)\text{-good}\}) \geq 1 - \varepsilon \qquad \text{if } d \geq d_0.$$

The quantity $\lambda/(1+\lambda)$, that appears in the statement of Lemma 2.1, is the equilibrium density of the independent flip process. We will use the notation

$$\bar{L}(\lambda) = L(\lambda/(1+\lambda)) = L(\lambda, \alpha) \qquad \text{for } \alpha \geq \lambda/(1+\lambda).$$

Note that in particular, for $L = \bar{L}(\lambda)$ and arbitrary $\delta > 0$, the configuration with all sites occupied in $H$ is $(L, \delta)$-good, if $d$ is large.

The proof of Lemma 2.1 will be given in Section 3. The next lemma, whose proof builds on Lemma 2.1, and will be postponed to Section 4, will allow us to use (2.1) to prove Proposition 2.1. Given $L$ and $\delta$, define the rescaling operator $\mathcal{R}_{L,\delta} : \{0,1\}^G \to \{0,1\}^{\mathbb{Z}^2}$ as follows:

$$(\mathcal{R}_{L,\delta}(\zeta))_i = \begin{cases} 1, & \text{if } \zeta \text{ is } (L, \delta)\text{-good on } H_i, \\ 0, & \text{otherwise,} \end{cases}$$

$i \in \mathbb{Z}^2$. Below $\succeq$ means stochastically larger than.

LEMMA 2.2 (Rescaling). *For any $\lambda > 0$ and $\varepsilon > 0$ there is $\delta > 0$ and finite $d_0$ such that for $L = \bar{L}(\lambda)$, $d \geq d_0$, and any $\zeta \in \{0,1\}^G$,*

$$\mathcal{R}_{L,\delta}(\eta^\zeta_{G;d^2}) \succeq T_{\varepsilon;1} \mathcal{R}_{L,\delta}(\zeta).$$

Iterating this lemma [using the Markov property of $(\eta^\zeta_{G;t})_{t \geq 0}$], and referring to the remark after Lemma 2.1, we learn that for any $\varepsilon > 0$ there is $\delta > 0$ such that for large enough $d$,

(2.2) $$\mathcal{R}_{L,\delta}(\eta^1_{G;nd^2}) \succeq T_{\varepsilon;n} \mathbf{1}, \qquad n = 1, 2, 3, \ldots.$$

Taking $\varepsilon < \varepsilon_c$ in (2.2) and combining this estimate with (2.1), proves that for any $\lambda > 0$, when $d$ is large enough, the process $(\eta^1_{G;t})_{t \geq 0}$ is $(L, \delta)$-good on $H$ at arbitrarily large times (multiples of $d^2$), provided that $L$ and $\delta$ have been properly chosen. Since the empty configuration on $H$ is clearly not $(L, \delta)$-good, when $d$ is large enough, we conclude then that $(\eta^1_{G;t})_{t \geq 0}$ survives. This shows that Proposition 2.1 will have been proved once we prove Lemmas 2.1 and 2.2.

**3. Almost product distributions concentrate on good configurations (proof of Lemma 2.1).** In order to prove Lemma 2.1, we introduce the following definition. A site $x \in H$ is said to be $(\zeta, K)$-dangerous if in the independent flip process $(\pi^\zeta_{H;t})_{t \geq 0}$ there is some $t \in [0, d^2]$ when $x$ has fewer than $K$ occupied neighbors. Define now

$$D^\zeta_K = \{x \in H : x \text{ is } (\zeta, K)\text{-dangerous}\}.$$

Note that the definition of $D^\zeta_K$ refers only to the independent flip process $(\pi^\zeta_{H;t})_{t \geq 0}$, and not to the threshold 2 contact process $(\eta^\zeta_{H;t})_{t \geq 0}$. The next lemma makes nevertheless the connection with that latter process.

LEMMA 3.1 (Only dangerous sites can cause trouble). *For any $\zeta \in \{0,1\}^H$, any $L > 0$, and any $K \geq L + 2$, on the event that $D^\zeta_K$ is $L$-thin, we have $\Delta^\zeta \subset D^\zeta_K$. In particular, if $D^\zeta_K$ is $L$-thin, then also $\Delta^\zeta$ is $L$-thin.*

PROOF. We argue by contradiction. (Informally: we consider the first violation of the claim, and show that it cannot be a violation of the claim.) Consider

$$s = \inf\{t \in [0, d^2] : \Delta_t^\zeta \not\subset D_K^\zeta\}$$

(inf $\varnothing = \infty$). If the lemma were false, we would have $s < d^2$ and there would be some $x \in H \setminus D_K^\zeta$ with the following properties. At $(x, s)$ there would be an up mark at which there would be a flip from 0 to 1 in the process $(\pi_{H;t}^\zeta)$, while in the process $(\eta_{H;t}^\zeta)$ the spin at $x$ would be 0 immediately before and immediately after time $s$. This implies that immediately before $s$ the site $x$ had fewer than 2 neighbors that were occupied in the process $(\eta_{H;t}^\zeta)$. But, by the definition of $s$, for $t < s$, any site $y$ that has $\eta_{H;t}^\zeta(y) \neq \pi_{H;t}^\zeta(y)$ must be $(\zeta, K)$-dangerous. Therefore, if we are on the event that $D_K^\zeta$ is $L$-thin, we learn that immediately before $s$ the site $x$ had at most $L$ neighbors where $\pi_{H;t}^\zeta$ differed from $\eta_{H;t}^\zeta$, and hence had fewer than $L + 2$ neighbors that were occupied in the process $(\pi_{H;t}^\zeta)$. This implies that $x \in D_K^\zeta$; a contradiction. □

Note that the next lemma refers only to the independent flip process $(\pi_{H;t}^\zeta)_{t\geq 0}$, and does not mention the threshold 2 contact process $(\eta_{H;t}^\zeta)_{t\geq 0}$. The arguments used to prove this lemma have some affinities with ideas of Balogh, Bollobas and Morris (2007).

LEMMA 3.2 (Controlling dangerous sites). *For any $\lambda > 0$ and any $\alpha \in (0, 1]$, there is $L = L(\lambda, \alpha) > 0$, such that for every $M > 0$, $K > 0$ and $\varepsilon > 0$, there is a finite $d_0' = d_0'(\lambda, \alpha, M, K, \varepsilon)$ such that if $\nu$ is an $M$-almost $\alpha$-product distribution over $\{0,1\}^H$, then*

$$\int d\nu(\zeta) \mathbb{P}(D_K^\zeta \text{ is } L\text{-thin}) \geq 1 - \varepsilon \qquad \text{if } d \geq d_0'.$$

PROOF. We will use the following notation for sites in $H = \{0,1\}^d$. For each $A \subset \{1, \ldots, d\}$, $e_A$ will denote the indicator of the set $A$, that is, the element of $H$ that is 1 on $A$ and 0 on $\{1, \ldots, d\} \setminus A$. Note that the neighborhood of the origin in $H$ is $\mathcal{N}_{H,0} = \{e_{\{1\}}, \ldots, e_{\{d\}}\}$.

By the symmetries of $H$, with no loss in generality,

$$\int d\nu(\zeta) \mathbb{P}(D_K^\zeta \text{ is not } L\text{-thin}) \leq 2^d \int d\nu(\zeta) \mathbb{P}(|D_K^\zeta \cap \mathcal{N}_{H,0}| \geq L)$$

(3.1)
$$\leq 2^d d^L \int d\nu(\zeta) \mathbb{P}(e_{\{1\}}, \ldots, e_{\{L\}} \in D_K^\zeta).$$

Since $\nu$ is an $M$-almost $\alpha$-product distribution,

$$\int d\nu(\zeta)\mathbb{P}(e_{\{1\}},\ldots,e_{\{L\}} \in D_K^\zeta) \leq \int d\mu_\alpha(\xi)\mathbb{P}(e_{\{1\}},\ldots,e_{\{L\}} \in D_{K+M}^\xi),$$
(3.2)

where $\mu_\alpha$ is a product distribution over $\{0,1\}^H$ with density $\alpha$.

For $i=1,\ldots,L$, consider now the set $N_i = \{e_{\{i,j\}} : j = L+1,\ldots,d\}$. Note that these $L$ sets are disjoint from each other, that they have cardinality $d-L$ and that for each $i$, $N_i \subset \mathcal{N}_{H,e_{\{i\}}}$. Let

$$T_i^{\xi,K+M} = \inf\left\{t \in [0,d^2]: \sum_{j=L+1}^d \pi_{H;t}^\xi(e_{\{i,j\}}) < K+M\right\},$$

$i=1,\ldots,L$. Since the event $\{e_{\{i\}} \in D_{K+M}^\xi\}$ is included in the event $\{T_i^{\xi,K+M} \leq d^2\}$, and the sets $N_i$ are disjoint from each other,

$$\int d\mu_\alpha(\xi)\mathbb{P}(e_{\{1\}},\ldots,e_{\{L\}} \in D_{K+M}^\xi) \leq \int d\mu_\alpha(\xi)\mathbb{P}\left(\bigcap_{i=1}^L \{T_i^{\xi,K+M} \leq d^2\}\right)$$
(3.3)
$$= \prod_{i=1}^L \int d\mu_\alpha(\xi)\mathbb{P}(T_i^{\xi,K+M} \leq d^2).$$

We want to estimate from above the probability that $T_i^{\xi,K+M} \leq d^2$. Set $W_i^\xi = \{t \in [0,d^2+1]: \sum_{j=L+1}^d \pi_{H;t}^\xi(e_{\{i,j\}}) < K+M\}$, and let $|W_i^\xi|$ be the Lebesgue measure of this set. Define also $T_i'$ as the instant of the first up mark in $N_i$ after $T_i^{\xi,K+M}$. Note that $T_i' - T_i^{\xi,K+M}$ is exponentially distributed with rate $(d-L)\lambda$, and is independent of $T_i^{\xi,K+M}$. Therefore, for $d \geq 1/\lambda$, using Markov's inequality,

$$\mathbb{E}(|W_i^\xi|) \geq \frac{1}{\lambda d}\mathbb{P}(T_i^{\xi,K+M} \leq d^2)\mathbb{P}\left(T_i' - T_i^{\xi,K+M} > \frac{1}{\lambda d}\right)$$
$$= \frac{1}{\lambda d}\ \mathbb{P}(T_i^{\xi,K+M} \leq d^2)\exp\left(-\frac{(d-L)\lambda}{d\lambda}\right)$$
$$\geq \frac{1}{e\lambda d}\mathbb{P}(T_i^{\xi,K+M} \leq d^2).$$

Hence, using Fubini's theorem,

$$\int d\mu_\alpha(\xi)\mathbb{P}(T_i^{\xi,K+M} \leq d^2)$$
(3.4)
$$\leq e\lambda d \int d\mu_\alpha(\xi)\mathbb{E}(|W_i^\xi|)$$
$$\leq e\lambda d(d^2+1)P\left(\text{Bin}\left(d-L,\alpha \wedge \frac{\lambda}{1+\lambda}\right) < K+M\right),$$

where $\text{Bin}(n,p)$ stands for a binomial random variable corresponding to $n$ independent attempts, each with probability $p$ of success. And we used the fact that in the independent flip process $(\pi_{H;t}^{\xi})$, each site is independently occupied at time $t$ with a probability that varies monotonically in $t$, starting from $\alpha$ at $t=0$, and converging to the equilibrium value $\lambda/(1+\lambda)$ as $t \to \infty$ (as well as the fact that binomial random variables are stochastically monotone increasing in the probability of success).

In order to choose $L(\lambda, \alpha)$, we recall the standard large deviation estimate

$$(3.5) \quad P\left(\text{Bin}\left(\frac{d}{2}, \alpha \wedge \frac{\lambda}{1+\lambda}\right) \leq \frac{d}{4}\left(\alpha \wedge \frac{\lambda}{1+\lambda}\right)\right) \leq \exp(-\gamma d),$$

where $\gamma = \gamma(\alpha \wedge \frac{\lambda}{1+\lambda}) > 0$. For a reason that soon will become clear, we take $L = L(\lambda, \alpha) = L(\lambda/(1+\lambda) \wedge \alpha)$ sufficiently large so that

$$(3.6) \quad \exp(-\gamma L) < \tfrac{1}{2}.$$

From (3.1)–(3.5) above,

$$\int d\nu(\zeta) \mathbb{P}(D_K^{\zeta} \text{ is not } L\text{-thin})$$

$$\leq 2^d d^L \left[e\lambda d(d^2+1)\right.$$

$$\left. \times P\left(\text{Bin}\left(d-L, \alpha \wedge \frac{\lambda}{1+\lambda}\right) < K + M\right)\right]^L$$

$$(3.7) \quad \leq 2^d d^L \left[e\lambda d(d^2+1)\right.$$

$$\left. \times P\left(\text{Bin}\left(\frac{d}{2}, \alpha \wedge \frac{\lambda}{1+\lambda}\right) \leq \frac{d}{4}\left(\alpha \wedge \frac{\lambda}{1+\lambda}\right)\right)\right]^L$$

$$\leq 2^d d^L [e\lambda d(d^2+1)\exp(-\gamma d)]^L$$

$$= (e\lambda d^2(d^2+1))^L (2\exp(-\gamma L))^d,$$

where in the second inequality, we suppose that $d \geq d_1$, for some large enough $d_1 = d_1(L(\lambda, \alpha), K+M)$.

Given $\varepsilon > 0$, it is clear that (3.6) and (3.7) imply the existence of a finite $d_0' = d_0'(\lambda, \alpha, M, K, \varepsilon)$ (larger than $d_1$ and $1/\lambda$), so that the claim stated in the lemma holds. $\square$

PROOF OF LEMMA 2.1. Take $L = L(\lambda, \alpha)$ from Lemma 3.2, and apply Lemma 3.2 with $K = L + 2$ and $\varepsilon$ replaced by $\varepsilon\delta$, to obtain the existence of $d_0 = d_0(\lambda, \alpha, M, \delta, \varepsilon)$ such that

$$\int d\nu(\zeta) \mathbb{P}(D_K^{\zeta} \text{ is not } L\text{-thin}) \leq \varepsilon\delta,$$

for $d \geq d_0$. From Lemma 3.1, we have for every $\zeta \in \{0,1\}^H$,

$$\mathbb{P}(\Delta^\zeta \text{ is not } L\text{-thin}) \leq \mathbb{P}(D_K^\zeta \text{ is not } L\text{-thin}).$$

On the set $\{\zeta : \zeta \text{ is not } (L,\delta)\text{-good}\}$, the probability on the left-hand side of the last display is larger than $\delta$, by the definition of $(L,\delta)$-good configurations. Therefore, integrating that inequality over $\nu$, and using the previous one, we obtain

$$\delta \nu(\{\zeta : \zeta \text{ is not } (L,\delta)\text{-good}\}) \leq \int d\nu(\zeta) \mathbb{P}(\Delta^\zeta \text{ is not } L\text{-thin}) \leq \varepsilon \delta.$$

Canceling $\delta$, we obtain the desired claim. $\square$

**4. Rescaling (proof of Lemma 2.2).** Set

$$\widehat{H} = H \cup H_{(0,1)} \cup H_{(1,0)}.$$

The main step in the proof of Lemma 2.2 will be the following lemma.

LEMMA 4.1 (Local rescaling estimate). *For any $\lambda > 0$ and $\varepsilon > 0$ there is $\delta > 0$ and finite $d_0$ such that for $d \geq d_0$ the following holds. If $\zeta \in \{0,1\}^{\widehat{H}}$, satisfies one of the following two conditions:*

(a) $\zeta$ *is* $(\bar{L}(\lambda), \delta)$*-good on* $H$,

*or*

(b) $\zeta$ *is* $(\bar{L}(\lambda), \delta)$*-good on* $H_{(0,1)}$ *and on* $H_{(1,0)}$,

*then*

$$\mathbb{P}(\eta^\zeta_{\widehat{H}; d^2} \text{ is } (\bar{L}(\lambda), \delta)\text{-good on } H) \geq 1 - \varepsilon.$$

Lemmas 4.3 and 4.4 below will contain the basic ingredients in the proof of Lemma 4.1.

LEMMA 4.2 (Goodness is monotone). *If $\zeta' \leq \zeta''$ are two configurations in $\{0,1\}^H$, then $\Delta^{\zeta''} \subset \Delta^{\zeta'}$. In particular, if for some $L > 0$ and $\delta > 0$, $\zeta'$ is $(L,\delta)$-good, then also $\zeta''$ is $(L,\delta)$-good.*

PROOF. If $x \in \Delta^{\zeta''}$, then there is some up mark at $x$ at some time $s \leq d^2$, with the following properties. At $(x,s)$ there is a flip from 0 to 1 in the process $(\pi^{\zeta''}_{H;t})$, while in the process $(\eta^{\zeta''}_{H;t})$ the spin at $x$ is 0 immediately before and immediately after time $s$. This implies that immediately before $s$ the site $x$ had fewer than 2 neighbors that were occupied in the process $(\eta^{\zeta''}_{H;t})$. Since, by attractivity, $\eta^{\zeta'}_{H;t} \leq \eta^{\zeta''}_{H;t}$, for all $t \geq 0$, it follows that also in the process $(\eta_{H;t}^{\zeta'})$ immediately before $s$ the site $x$ was vacant and had fewer than 2 neighbors that were occupied. Therefore, $x$ will still be vacant in the process $(\eta_{H;t}^{\zeta'})$, immediately after $s$, while in the process $(\pi_{H;t}^{\zeta'})$ it will be occupied. This implies that $x \in \Delta^{\zeta'}$. □

LEMMA 4.3 (Goodness is likely to be preserved). *For any $\lambda > 0$, $M > 0$, and $\delta > 0$, for each $\zeta \in \{0,1\}^H$ that is $(M,\delta)$-good and $t \in [d^2/4, d^2]$, we have:*

(i) *The total variation distance between the law of $\eta_{H;t}^{\zeta}$ and a certain $M$-almost $\lambda/(1+\lambda)$-product distribution is bounded above by $2\delta$, for large $d$.*

(ii) *For any $\delta' > 0$,*

$$\mathbb{P}(\eta_{H;t}^{\zeta} \text{ is } (\bar{L}(\lambda), \delta')\text{-good}) \geq 1 - 4\delta \qquad \text{for large } d.$$

PROOF. Let $\xi$ be a random configuration in $\{0,1\}^H$, corresponding to a product distribution with density $\lambda/(1+\lambda)$. This is the equilibrium distribution for the independent flip process $(\pi_{H,t})$, so that $\xi_t = \pi_{H,t}^{\xi}$ has the same $\lambda/(1+\lambda)$-product distribution for each $t \geq 0$. Define now

$$\xi_t^{M,\zeta}(x) = \begin{cases} \xi_t(x), & \text{if } x \in H \setminus \Delta_t^{\zeta} \text{ or } \Delta_t^{\zeta} \text{ is not } M\text{-thin,} \\ 0, & \text{if } x \in \Delta_t^{\zeta} \text{ and } \Delta_t^{\zeta} \text{ is } M\text{-thin.} \end{cases}$$

Note that, for arbitrary $t \geq 0$, $\xi_t^{M,\zeta}$ has an $M$-almost $\lambda/(1+\lambda)$-product distribution.

The event $A$ that there is at least one Poisson mark, either up or down, at each site $x \in H$ during the time interval $[0, d^2/4]$ has

$$\mathbb{P}(A^c) \leq 2^d \exp(-d^2/4) \leq \delta,$$

for large $d$. And on $A$, $\pi_{H;t}^{\zeta} = \xi_t$, for $t \geq d^2/4$.

Let $B$ denote the event that $\Delta^{\zeta}$ is $M$-thin. From the assumption on $\zeta$, we have

$$\mathbb{P}(B^c) \leq \delta.$$

On $AB$ we have

$$\xi_t^{M,\zeta}(x) = \begin{cases} \xi_t(x) = \pi_{H;t}^{\zeta}(x) = \eta_{H;t}^{\zeta}(x), & \text{if } x \in H \setminus \Delta_t^{\zeta} \\ 0 = \eta_{H;t}^{\zeta}(x), & \text{if } x \in \Delta_t^{\zeta} \end{cases} \Bigg\} = \eta_{H;t}^{\zeta}(x).$$

This means that $\eta_{H;t}^{\zeta} = \xi_t^{M,\zeta}$, on the event $AB$, which has $\mathbb{P}((AB)^c) \leq 2\delta$, when $d$ is large, completing the proof of (i).

The claim (ii) follows from (i) and Lemma 2.1, by possibly taking $d$ larger. □

LEMMA 4.4 (Infecting with goodness). *For any $\lambda > 0$, $M > 0$, and $\beta > 0$, if $\psi$ is a random configuration on $\{0,1\}^{\widehat{H}}$ that has restrictions to $H_{(0,1)}$ and $H_{(1,0)}$ that are independent and each has an $M$-almost $\beta$-product distribution, then we have:*

  (i) *$\eta^{\psi}_{\widehat{H};1}$ has a restriction to $H$ that is stochastically larger than an $8M$-almost $\alpha$-product distribution, where $\alpha = \alpha(\lambda, \beta) > 0$.*

  (ii) *For any $\varepsilon > 0$ and $\delta > 0$,*

$$\mathbb{P}(\eta^{\psi}_{\widehat{H};1} \text{ is } (L(\lambda, \alpha), \delta)\text{-good on } H) \geq 1 - \varepsilon \qquad \text{for large } d.$$

PROOF. For each $i = (i_1, i_2) \in \mathbb{Z}^2$, partition $H_i$ into the sets

$$H_{i,y} = \{2i_1, 2i_1 + 1\} \times \{2i_2, 2i_2 + 1\} \times \{y\}, \qquad y \in \{0,1\}^{d-2}.$$

From the hypothesis, we have, for $i = (0,1)$ and $i = (1,0)$, a set $S_i \subset H_i$ that is $M$-thin in $H_i$, and such that the restriction of $\psi$ to $H_i$ is identical on $H_i \setminus S_i$ to a random configuration $\xi$ on $H_{(0,1)} \cup H_{(1,0)}$ that has a $\beta$-product distribution.

Define

$$I = \{y \in \{0,1\}^{d-2} : H_{(0,1),y} \cap S_{(0,1)} \neq \varnothing \text{ or } H_{(1,0),y} \cap S_{(1,0)} \neq \varnothing\},$$

so that for $y \in \{0,1\}^{d-2} \setminus I$, $\psi$ is identical to $\xi$ on $H_{(0,1),y} \cup H_{(1,0),y}$.

For $y = (y_3, \ldots, y_d) \in \{0,1\}^{d-2}$, consider the following event $A_y$: $\xi(x) = 1$ for $x \in H_{(0,1),y} \cup H_{(1,0),y}$, between times $0$ and $1$, there is no down mark in $H_{(0,0),y} \cup H_{(0,1),y} \cup H_{(1,0),y}$ and between times $0$ and $1$ there is an up mark at each one of the sites $(1, 1, y_3, \ldots, y_d)$, $(1, 0, y_3, \ldots, y_d)$, $(0, 1, y_3, \ldots, y_d)$ and $(0, 0, y_3, \ldots, y_d)$, in that order, or in the order $(1, 1, y_3, \ldots, y_d)$, $(0, 1, y_3, \ldots, y_d)$, $(1, 0, y_3, \ldots, y_d)$, $(0, 0, y_3, \ldots, y_d)$. Then the events $A_y$, $y \in \{0,1\}^{d-2}$, are i.i.d., with

$$\mathbb{P}(A_y) \geq \beta^8 e^{-12}(1 - e^{-\lambda/4})^4 = \alpha'(\lambda, \beta) > 0.$$

Define also a random configuration $\phi$ on $H$, by setting

$$\phi(x) = \begin{cases} 1, & \text{if } x \in H_{(0,0),y} \text{ for some } y \text{ for which } A_y \text{ occurs,} \\ 0, & \text{otherwise.} \end{cases}$$

Set

$$S = \bigcup_{y \in I} H_{(0,0),y}.$$

We have now the following estimate, which is central to our argument, and follows immediately from the definitions above, and the way $(\eta^{\psi}_{\widehat{H};t})$ is constructed from the Poisson marks:

$$\eta^{\psi}_{\widehat{H};1}(x) \geq \phi(x) \qquad \text{if } x \in H \setminus S.$$

Claim (i) in the lemma will therefore follow, once we argue that the law of $\phi$ is stochastically larger than a product distribution with some density $\alpha > 0$, and that $S$ is $8M$-thin.

For the first of these tasks, it is enough to take $\alpha = \alpha(\lambda, \beta) = \alpha'((\lambda, \beta)/4$. Indeed, a product distribution over $H$ with density $\alpha$ gives, for each $y \in \{0,1\}^{d-2}$, probability $(1 - \alpha'/4)^4 \geq 1 - \alpha'$ to the event that all 4 sites in $H_{(0,0),y}$ are vacant. And this is enough, since in $\phi$, when $H_{(0,0),y}$ is not fully vacant, it is fully occupied.

We now argue that $S$ is an $8M$-thin subset of $H$. For this, let $S'_{(0,1)}$ and $S'_{(1,0)}$ be, respectively, the translates of $S_{(0,1)}$ and $S_{(1,0)}$ to $H$. Then $S' = S'_{(0,1)} \cup S'_{(1,0)}$ is $2M$-thin. To obtain $S$ from $S'$, we introduce the mappings $\sigma_k : H \to H$, $k = (k_1, k_2) \in \{0,1\}^2$, given by

$$(\sigma_k(x))_j = \begin{cases} 1 - x_j, & \text{if } j \in \{1,2\} \text{ and } k_j = 1, \\ x_j, & \text{if } j \in \{3, \ldots, d\} \text{ or } k_j = 0. \end{cases}$$

These 4 mappings are automorphisms of $H$ (as a graph), and

$$S = \bigcup_{k \in \{0,1\}^2} \sigma_k(S').$$

So $S$ is the union of 4 sets that are $2M$-thin, and hence it is $8M$-thin.

This completes the proof of claim (i).

Claim (ii) is now a consequence of (i) and Lemmas 2.1 and 4.2. □

PROOF OF LEMMA 4.1. In case condition (a) holds, the claim in the lemma is immediate from (ii) in Lemma 4.3 [with $M = \bar{L}(\lambda)$, $\delta' = \delta$], provided that we take $\delta \leq \varepsilon/4$.

In case condition (b) holds, we will show next that the claim in the lemma holds, provided that we take $\delta \leq \varepsilon/12$. For this purpose, we will compare the process $(\eta^\zeta_{\widehat{H};t})_{t \in [0,d^2]}$, with the process $(\check{\eta}^\zeta_{\widehat{H};t})_{t \in [0,d^2]}$ in which interaction is only allowed among the hypercubes $H$, $H_{(0,1)}$ and $H_{(1,0)}$ during the time interval $[d^2/2 - 1, d^2/2]$.

The random configuration $\check{\eta}^\zeta_{\widehat{H};d^2/2-1}$ has restrictions to $H$, $H_{(0,1)}$ and $H_{(1,0)}$ that are independent. From (i) in Lemma 4.3 [in its translated versions for $H_{(0,1)}$ and $H_{(1,0)}$, with $M = \bar{L}(\lambda)$], we learn that, for large $d$, this random configuration has restrictions to $H_{(0,1)}$ and $H_{(1,0)}$ each with a law that is at total variation distance at most $2\delta$ from a $\bar{L}(\lambda)$-almost $\lambda/(1+\lambda)$-product distribution. Therefore the random configuration $\check{\eta}^\zeta_{\widehat{H};d^2/2-1}$ has a law that is at total variation distance at most $4\delta$ from a random configuration $\psi$ on $\widehat{H}$ that satisfies the conditions of Lemma 4.4, with $M = \bar{L}(\lambda)$ and $\beta = \lambda/(1+\lambda)$. [See, e.g., Lemma (6.2), Section 2.6, page 139 of Durrett (1996), for the fact that the total variation distance is added when measures are multiplied.]

THRESHOLD CONTACT PROCESSES 15

Using now the Markov property of $(\check{\eta}^{\zeta}_{\widehat{H};t})_{t\in[0,d^2]}$ (applied at time $d^2/2-1$), and (ii) in Lemma 4.4 (with the time interval $[0,1]$ translated to $[d^2/2-1, d^2/2]$), we learn that for an appropriate $\alpha > 0$,

$$\mathbb{P}(\check{\eta}^{\zeta}_{\widehat{H};d^2/2} \text{ is } (L(\lambda,\alpha),\delta)\text{-good on } H)$$
$$\geq \mathbb{P}(\eta^{\psi}_{\widehat{H};1} \text{ is } (L(\lambda,\alpha),\delta)\text{-good on } H) - 4\delta \geq 1 - \varepsilon/3 - 4\delta,$$

for large enough $d$.

Using once more the Markov property of $(\check{\eta}^{\zeta}_{\widehat{H};t})_{t\in[0,d^2]}$ (this time applied at time $d^2/2$), and using (ii) in Lemma 4.3 [this time with $M = L(\lambda,\alpha)$, $t = d^2/2$ and $\delta' = \delta$], we learn now that

$$\mathbb{P}(\check{\eta}^{\zeta}_{\widehat{H};d^2} \text{ is not } (\bar{L}(\lambda),\delta)\text{-good on } H)$$
$$\leq \mathbb{P}(\check{\eta}^{\zeta}_{\widehat{H};d^2/2} \text{ is not } (L(\lambda,\alpha),\delta)\text{-good on } H)$$
$$+ \mathbb{P}(\check{\eta}^{\zeta}_{\widehat{H};d^2} \text{ is not } (\bar{L}(\lambda),\delta)\text{-good on } H | \check{\eta}^{\zeta}_{\widehat{H};d^2/2} \text{ is } (L(\lambda,\alpha),\delta)\text{-good on } H)$$
$$\leq \varepsilon/3 + 4\delta + 4\delta \leq \varepsilon,$$

for large enough $d$. The claim in the lemma now follows by using attractivity and Lemma 4.2, when comparing the process $(\eta^{\zeta}_{\widehat{H};t})_{t\in[0,d^2]}$ with the process $(\check{\eta}^{\zeta}_{\widehat{H};t})_{t\in[0,d^2]}$. $\square$

PROOF OF LEMMA 2.2. We will combine Lemma 4.1, with Theorem 0.0(i) of Liggett, Schonmann and Stacey (1997). For this purpose, for $i \in \mathbb{Z}^2$, set

$$\widehat{H}_i = H_i \cup H_{i+(0,1)} \cup H_{i+(1,0)}$$

and define the random variables $\widetilde{X}^{\zeta}_i$, $i \in \mathbb{Z}^2$, by setting $\widetilde{X}^{\zeta}_i = 1$ in case $(\mathcal{R}_{L,\delta}(\eta^{\zeta}_{\widehat{H}_i;d^2}))_i = 1$, or $(\mathcal{R}_{L,\delta}(\zeta))_i = (\mathcal{R}_{L,\delta}(\zeta))_{i+(0,1)}(\mathcal{R}_{L,\delta}(\zeta))_{i+(1,0)} = 0$; and $\widetilde{X}^{\zeta}_i = 0$, otherwise. Then $(\widetilde{X}^{\zeta}_i)_{i\in\mathbb{Z}^2}$ has a finite range of dependence. Therefore, Theorem 0.0(i) of Liggett, Schonmann and Stacey (1997) assures us that given $\varepsilon > 0$, there exists $\tilde{\varepsilon} > 0$ such that if

$$\inf_{i \in \mathbb{Z}^2} \mathbb{P}(\widetilde{X}^{\zeta}_i = 1) \geq 1 - \tilde{\varepsilon}, \tag{4.1}$$

then there is an i.i.d. random field $(X_i)_{i\in\mathbb{Z}^2}$, with $P(X_i = 0) = \varepsilon$, $P(X_i = 1) = 1 - \varepsilon$, that lies stochastically below $(\widetilde{X}^{\zeta}_i)_{i\in\mathbb{Z}^2}$.

On the other hand, Lemma 4.1 (with $\tilde{\varepsilon}$ in place of the $\varepsilon$ there), assures us that (4.1) is indeed satisfied, provided that we take $\delta > 0$ properly, and $d$ large enough.

Define the following random transformation, $\widetilde{T}$, on $\{0,1\}^{\mathbb{Z}^2}$:

$$(\widetilde{T}\Psi)_i = \begin{cases} \widetilde{X}_i^{\zeta}, & \text{if } \Psi_i = 1 \text{ or } \Psi_{i+(0,1)} = \Psi_{i+(1,0)} = 1, \\ 0, & \text{otherwise,} \end{cases}$$

for $\Psi \in \{0,1\}^{\mathbb{Z}^2}$, $i \in \mathbb{Z}^2$. Under the conditions above on $\delta$ and $d$, we have now,

$$\mathcal{R}_{L,\delta}(\eta_{G;d^2}^{\zeta}) \succeq \widetilde{T}\mathcal{R}_{L,\delta}(\zeta) \succeq T_{\varepsilon;1}\mathcal{R}_{L,\delta}(\zeta).$$

This completes the proof. □

**5. Extensions, related results, notes.** We will use the following standard notation, for the elements of the canonical base of $\mathbb{Z}^d$: $e_1 = (1,0,0,\ldots,0,0)$, $e_2 = (0,1,0,\ldots,0,0)$, ..., $e_d = (0,0,0,\ldots,0,1)$. For $x \in \mathbb{Z}^d$, we say that its nearest neighbors $y$ and $z$ are colinear with $x$ if $\{y,z\} = \{x - e_j, x + e_j\}$ for some $j \in \{1,\ldots,d\}$.

We will be considering other processes on $\mathbb{Z}^d$ which can be constructed with the same system of Poisson marks used for the threshold contact processes. As usual, we say that model A dominates model B, if when they are both constructed in this way and are started from the same configuration, model A is at any time above model B.

(1) First we explain how the proof of Theorem 1.1 should be modified for arbitrary $\theta$. One uses the following generalization of (2.1). Let $X_{i,n}$, $i \in \mathbb{Z}^d$, $n \in \mathbb{Z}$, be a collection of i.i.d. random variables, with $P(X_{i,n} = 0) = \varepsilon$, $P(X_{i,n} = 1) = 1 - \varepsilon$, where $\varepsilon \in [0,1]$ is a noise parameter. Define now the following family of random transformations, $T_{\varepsilon;n}$, on $\{0,1\}^{\mathbb{Z}^d}$, indexed by $n \in \{0,1,2,3,\ldots\}$. $T_{\varepsilon;0}$ is the identity and recursively, for $n = 1, 2, \ldots$,

$$(T_{\varepsilon;n}\Psi)_i = \begin{cases} X_{i,n}, & \text{if } (T_{\varepsilon;n-1}\Psi)_i = 1, \text{ or } (T_{\varepsilon;n-1}\Psi)_{i+e_j} = 1, \ j = 1,\ldots,d, \\ 0, & \text{otherwise,} \end{cases}$$

for $\Psi \in \{0,1\}^{\mathbb{Z}^d}$, $i \in \mathbb{Z}^d$. Denote by **1** the element of $\{0,1\}^{\mathbb{Z}^d}$ identically 1. Toom (1974) proved the following result [see also Bramson and Gray (1991), for an alternative proof]: There exists $\varepsilon_c > 0$ such that

(5.1) $\quad$ for $\varepsilon < \varepsilon_c$, $\quad \liminf_{n \to \infty} P((T_{\varepsilon;n}\mathbf{1})_0 = 1) > 0$.

For our purposes, we set

$$G = \mathbb{Z}^\theta \times \{0,1\}^{d-\theta} = \bigcup_{i \in \mathbb{Z}^\theta} H_i,$$

where

$$H_i = \{2i_1, 2i_1 + 1\} \times \cdots \times \{2i_\theta, 2i_\theta + 1\} \times \{0,1\}^{d-\theta}, \qquad i = (i_1, \ldots, i_\theta) \in \mathbb{Z}^\theta.$$

One then proves the analogue of Proposition 2.1, exactly in the same way as that proposition was proved, except for obvious changes [e.g., (5.1) is used in its $\theta$-dimensional version, in Lemma 3.1 the condition reads $K \geq L + \theta$, $\widehat{H} = H \cup (\bigcup_{j=1,\ldots,\theta} H_{e_j})$, in (i) in Lemma 4.4, the quantity $8M$ becomes $\theta 2^\theta M$, etc.].

(2) Even when the threshold $\theta$ contact process survives started from all sites occupied, it may happen that for small $p > 0$, $\eta_t^p \Rightarrow \delta_0$, as $t \to \infty$. By attractivity, if this happens for some value of $p$, it will also happen for smaller values of $p$. This leads to the following definitions:

$$p_c(\lambda) = \sup\{p \in [0,1] : \eta_t^p \Rightarrow \delta_0, \text{ as } t \to \infty\}$$
$$= \sup\{p \in [0,1] : \rho_t^p \to 0 \text{ as } t \to \infty\},$$
$$\bar{\lambda} = \inf\{\lambda : p_c(\lambda) = 0\} = \inf\left\{\lambda : \text{for all } p > 0, \inf_{t \geq 0} \rho_t^p > 0\right\}.$$

In case $\theta = 1$, one can use duality to easily show that whenever the process survives, it survives starting from any product distribution with positive density. In particular, $\lambda_c = \bar{\lambda}$.

Chen (1994) proved a result that implies that if $\theta = 2$, $d \geq 3$, then

(5.2) $$\bar{\lambda} < \infty.$$

Specifically, Chen (1994) studied a spin flip model on $\mathbb{Z}^d$ that he referred to as the system with symmetric sexual reproduction. That model is the $\theta = 2$ case of the model that we call the modified threshold $\theta$ contact process, and is defined by the following rates. At $x \in \mathbb{Z}^d$, at time $t \geq 0$:

- 1 flips to 0 at rate 1.
- 0 flips to 1 at rate $\lambda$ in case there are at least $\theta$ nearest neighbors of $x$, no two of which are colinear with $x$, in state 1 at time $t$; and at rate 0 otherwise.

Theorem 4 of Chen (1994) shows that when $d \geq 3$, the modified threshold 2 contact process on $\mathbb{Z}^d$, with a large enough $\lambda$, survives for any initial density $p > 0$. By domination, the same therefore holds for the threshold 2 contact process on $\mathbb{Z}^d$, that is, (5.2) holds.

In Chen (1992) the following result, that contrasts with (5.2), was proved. In $d = 2$ the modified threshold 2 contact process with arbitrarily large $\lambda$, started from a product measure with small enough density $p > 0$ (how small depends on $\lambda$) dies out. This means that for this model $\bar{\lambda} = \infty$. Presumably the same holds for the threshold 2 contact process in dimension 2.

Interestingly enough, the qualitative behavior of the threshold 2 contact process in high dimensions is distinct from the behavior of the associated mean-field model, analyzed in Fontes and Schonmann (2008). In that mean-field model, $p_c(\lambda) > 0$ for all $\lambda > 0$, with $\lim_{\lambda \to \infty} p_c(\lambda) = 0$. In particular, $\bar{\lambda} =$ $\infty$. Also the behavior of the threshold 2 contact process on a homogeneous tree is distinct from that of both the process on $\mathbb{Z}^d$ with large $d$ and the mean-field model. For the models on the homogeneous trees with $\theta \geq 2$, Fontes and Schonmann (2008) observed that even $\lim_{\lambda \to \infty} p_c(\lambda) > 0$.

In connection with the facts reviewed above it seems very interesting to decide whether, for $\theta = 2$, in $d \geq 3$ the critical points $\lambda_c$ and $\bar{\lambda}$ are distinct from each other.

(3) It is worth pointing out that Theorem 1.1 extends also to the modified threshold $\theta$ contact process, as defined in (2) above. This is so simply because the modified threshold $\theta$ contact process dominates the threshold $2\theta - 1$ contact process.

(4) The following strengthening of Theorem 1.1 can be proved with very little additional work.

THEOREM 5.1. *For any $\theta$, and any $\lambda > 0$,*

$$\lim_{d \to \infty} p_c(\lambda) = 0.$$

To prove this theorem, first we prove that for arbitrary $p > 0$, $\varepsilon' > 0$ and $\delta' > 0$, for $L = \bar{L}(\lambda)$,

(5.3) $$\mathcal{R}_{L,\delta'}(\eta^p_{G;d^2}) \succeq T_{\varepsilon';1}\mathbf{1}$$

for large $d$. This claim is a direct consequence of Lemma 2.1 [with $\alpha = p$, $M$ irrelevant, $\varepsilon = \varepsilon'/2$, $L = L(\lambda, p)$ and $\delta = \varepsilon'/8$], combined with (ii) in Lemma 4.3 [with $M = L(\lambda, p)$, $\delta = \varepsilon'/8$ and $t = d^2$], and with the fact that the right-hand side of (5.3) has a product distribution with density $\varepsilon'$.

Combining (5.3) with $n - 1$ iterates of Lemma 2.2 [using the Markov property of $(\eta^p_{G;t})_{t \geq 0}$], we learn that for any $\varepsilon > 0$ there is $\delta > 0$ such that for $L = \bar{L}(\lambda)$ and large enough $d$,

(5.4) $$\mathcal{R}_{L,\delta}(\eta^p_{G;nd^2}) \succeq T_{\varepsilon;n}\mathbf{1}, \qquad n = 1, 2, 3, \ldots.$$

Theorem 5.1 now follows from (5.4) in the same way as Theorem 1.1 followed from (2.2).

(5) We end this section with a proof that must be well known, but does not seem to be in print. We observe how (5.1) can be used to prove that for every $d \geq 2$, the asymmetric model defined by the rates below has $\lambda_c < \infty$. At $x \in \mathbb{Z}^d$, at time $t \geq 0$:

- 1 flips to 0 at rate 1.
- 0 flips to 1 at rate $\lambda$ in case $x + e_j$, $j = 1, \ldots, d$, are all in state 1 at time $t$; and at rate 0 otherwise.

We will use $(\vec{\eta}_t)_{t\geq 0}$ to denote this process. Define the random variables $X_{i,n}$, $i \in \mathbb{Z}^d$, $n = 1, 2, \ldots$, by $X_{i,n} = 1$ if at $i$, during the time interval $[(n-1)/\sqrt{\lambda}, n/\sqrt{\lambda}]$ there is no down mark and there is an up mark; $X_{i,n} = 0$, otherwise. Then the $X_{i,n}$ are i.i.d., and $\varepsilon = \mathbb{P}(X_{i,n} = 0)$ can be made arbitrarily small, by taking $\lambda$ large. It is also clear that

$$\vec{\eta}^{\mathbf{1}}_{n/\sqrt{\lambda}} \succeq T_{\varepsilon;n}\mathbf{1}, \qquad n = 1, 2, 3, \ldots.$$

Therefore (5.1) implies survival of $(\vec{\eta}_t)_{t\geq 0}$ when $\lambda$ is large enough. By domination, we learn also that for the threshold contact processes on $\mathbb{Z}^d$, with $\theta \leq d$, $\lambda_c < \infty$, as mentioned in the Introduction.

**Acknowledgment.** Roberto H. Schonmann is pleased to thank the warm hospitality of the Institut de Mathématiques, École Polytechnique Fédérale de Lausanne, where this work was started.

Institut de Mathématiques
École Polytechnique Fédérale
Station 8, 1015 Lausanne
Switzerland
E-mail: thomas.mountford@epfl.ch

Department of Mathematics
University of California at
  Los Angeles
Los Angeles, California 90095
USA
E-mail: rhs@math.ucla.edu